\documentclass[10pt]{amsart}
\usepackage{amsmath, amsfonts, amsthm, amssymb, verbatim,dsfont}
\usepackage{graphicx}
\usepackage[mathcal]{euscript}
\usepackage{amstext}
\usepackage{dsfont} 
\usepackage{amssymb,mathrsfs}
\usepackage[usenames]{color}
\theoremstyle{plain}
        \newtheorem{theorem}{Theorem}[section]
        
        \newtheorem{lemma}[theorem]{Lemma}
        
\theoremstyle{definition}

\theoremstyle{plain}
                
\numberwithin{equation}{section}

\usepackage{amsmath}
\usepackage{verbatim}
\newcommand \be           {\begin{equation}}
\newcommand \ee            {\end{equation}}

\newcommand \Fcal           {\mathcal{F}}
\newcommand \RR           {\mathbb{R}}
\newcommand \NN           {\mathbb{N}}

\newcommand \CC           {\mathbb{C}}
\newcommand \BB           {\mathbb{B}}

\newcommand \Pbold           {\mathbf{P}} 
\newcommand \PP \Pbold

\newcommand \del           \partial
\newcommand \eps            \epsilon

\newcommand \loc        {{\mathrm{loc}}}

\DeclareMathOperator    \Id {Id}

\newcommand \ws     {\mathrel{\mathop{\rightharpoonup}\limits^{*}}}
\definecolor{gray}{gray}{0.4}
\newcommand{\unj}{u^n_j}
\newcommand{\unjj}{u^n_{j+1}}
\newcommand{\unjm}{u^n_{j-1}}
\newcommand{\unnj}{u^{n+1}_j}

\newcommand{\vnj}{v^n_j}
\newcommand{\vnjj}{v^n_{j+1}}
\newcommand{\vnjm}{v^n_{j-1}}
\newcommand{\vnnj}{v^{n+1}_j}

\newcommand{\wnj}{w^n_j}
\newcommand{\wnjj}{w^n_{j+1}}
\newcommand{\wnjm}{w^n_{j-1}}
\newcommand{\wnnj}{w^{n+1}_j}

\def\XXint#1#2#3{{\setbox0=\hbox{$#1{#2#3}{\int}$}
\vcenter{\hbox{$#2#3$}}\kern-.5\wd0}}

\def\build#1_#2^#3{\mathrel{
\mathop{\kern 0pt#1}\limits_{#2}^{#3}}}

\begin{document}

\title[short wave long wave interactions in a viscoelastic medium]{A nonlinear model describing a short wave long wave interaction in a viscoelastic medium}

\author[P. Amorim \and J.P. Dias]{Paulo Amorim$^1$ \and Jo\~ ao Paulo Dias$^1$}

\date{\today}

\footnotetext[1]{Centro de Matem\'atica e Aplica\c c\~oes
Fundamentais, FCUL, Av. Prof. Gama Pinto 2,
1649-003 Lisboa, Portugal. E-mail: {\tt pamorim@ptmat.fc.ul.pt, dias@ptmat.fc.ul.pt}}

\begin{abstract} 
In this paper we introduce a system coupling a nonlinear Sch\-r\"o\-din\-ger equation with a system of viscoelasticity, modeling the interaction between short and long waves, acting for instance on media like plasmas or polymers. We prove the existence and uniqueness of local (in time) strong solutions and the existence of global weak solutions for the corresponding Cauchy problem. In particular we extend previous results in [Nohel \emph{et. al.}, Commun. Part. Diff. Eq., 13 (1988)] for the quasilinear system of viscoelasticity. We finish with some numerical computations to illustrate our results.
\end{abstract}

 
\maketitle

\section{Introduction}
In \cite{Benney}, D.J. Benney initiated 
the study of mathematical models describing the interaction between short and long waves in fluids, namely capillary and gravity waves or internal and surface waves. This has been developed by many authors, such as M. Tsutsumi and S. Hatano in the pioneering papers \cite{TH1,TH2}, and more recently, in the framework of quasilinear systems, in \cite{DF1,DF2,DFF,DFr,DFO1} and \cite{DFO2}. Their numerical study was initiated in \cite{AF1} and developed in \cite{AF2}.

In this paper, we will apply the ideas of these authors to the case of viscoelastic fluids, arising, for instance, in plasma physics (in the study of helioseismology, cf. \cite{JCD}), or in the study of polymers \cite{LP}. To this purpose, we introduce a model coupling the nonlinear Sch\-r\"o\-din\-ger equation (modeling the short waves) with a quasilinear system of viscoelasticity (modeling the long waves). This last system was studied in \cite{NRT}. Thus, we analyze the following Cauchy problem,
\be
\label{10}
\left\{
\aligned
&i u_t + u_{xx} = uv + \alpha |u|^2 u
\\
&v_t = w_x
\\
& w_t = (\sigma(v))_x + (|u|^2)_x + \int_0^t k(t- \tau) \big[ (\sigma(v))_x + (|u|^2)_x \big] (\cdot,\tau) \,d\tau,
\endaligned
\right.
\ee 
where $x \in \RR$, $t\ge 0$, $\alpha$ is a real constant, and $i = \sqrt{-1}$. Here, $u(x,t) \in \CC$ is the envelope of the short waves, $v(x,t)\in \RR$ is the deformation gradient, $w(x,t) \in \RR$ the velocity of the long waves,  $\sigma \in C^3(\RR)$ is the stress function verifying $\sigma' \ge \sigma_0 >0$ (hyperbolicity), and $k$ is a given $C^1$ kernel. We consider \eqref{10} with the initial data
\be
\label{20}
\aligned
u_0,\, v_0,\, w_0 \in H^1(\RR).
\endaligned
\ee
The integral term in the right-hand side of \eqref{10} represents the memory effects due to the viscoelastic structure of the fluid.

We now recall a transformation due to R.C. MacCamy \cite{MacC} (see also \cite{NRT}): Let $q(t)$ be the resolvent kernel associated with $k$, i.e., $q$ is the solution of the linear Volterra equation
\be
\label{30}
\aligned
q(t) + \int_0^t k(t - \tau) q(\tau) \,d\tau, \qquad \tau\ge 0.
\endaligned
\ee
Convolving the third equation in \eqref{10} with $q(t)$, it is not difficult to obtain for smooth solutions
\be
\label{40}
\aligned
&\int_0^t k(t - \tau) \big[ (\sigma(v))_x + (|u|^2)_x \big] (x,\tau) \,d\tau  
 = \int_0^t q(t - \tau) w_t(x, \tau) \,d\tau 
\\
&\qquad= q(0) w(x,t) - q(t) w_0(x) + \int_0^t q'(t-\tau) w(x,\tau) \,d\tau.
\endaligned
\ee
Thus, for smooth solutions, \eqref{10}\eqref{20} is equivalent to the Cauchy problem
\be
\label{50}
\left\{
\aligned
&i u_t + u_{xx} = uv + \alpha |u|^2 u
\\
&v_t = w_x && x\in\RR, \quad t\ge0,
\\
& w_t = (\sigma(v))_x + (|u|^2)_x + \Fcal(w),
\endaligned
\right.
\ee 
\be
\label{60}
\aligned
u(x,0) = u_0 (x), \quad v(x,0) = v_0 (x), \quad w(x,0) = w_0(x),
\endaligned
\ee
where
\be
\label{70}
\Fcal(w)(x,t) = q(0) w(x,t) - q(t) w_0(x) + \int_0^t q'(t-\tau) w(x,\tau) \,d\tau.
\ee

In the case where $u\equiv 0$, that is, the nonlinear viscoelasticity system, the existence of a global (in time) weak solution for initial data in $L^\infty (\RR) \cap L^2 (\RR)$ and $k \in C^1([0, +\infty))$ was proved in \cite{NRT} (see also \cite{CD} for a different model), by the vanishing viscosity method applied to both variables $v$ and $w$ and the compensated compactness method \cite{T,DP} based on ideas introduced by C.~Dafermos in \cite{D1} (see also \cite{DF2}, and for related results \cite{D3,DF1}).

In this paper, we start by proving in Section \ref{S20} the existence and uniqueness for local (in time) strong solutions for \eqref{50},\eqref{60} with initial data $(u_0, v_0, w_0) \in H^3(\RR) \times (H^2(\RR))^2$ by applying a variant of Kato's theorem \cite{K}. In Section~\ref{S30} and for a special class of stress functions $\sigma$, we apply the physical vanishing viscosity method (that is, only in the velocity variable $w$) and a variant of the compensated compactness method introduced by D. Serre and J. Shearer in \cite{SS} (see also \cite{DFF}) to obtain the existence of a global (in time) weak solution for \eqref{50},\eqref{60} with initial data $(u_0, v_0, w_0) \in (H^1(\RR))^3.$ Finally, in Section~\ref{S40}, we present some numerical simulations to illustrate the behavior of the solutions in a special case.

\section{Local (in time) existence of strong solutions}
\label{S20}
Let $u_0 \in H^3(\RR)$, $v_0 \in H^2(\RR)$, $w_0 \in H^2(\RR)$. To study the Cauchy problem \eqref{50},\eqref{60}, we introduce the Riemann invariants
\[
\aligned
l = w + \int_0^v \sqrt{\sigma'(\xi)} \,d\xi, \qquad r = w - \int_0^v \sqrt{\sigma'(\xi)} \,d\xi.
\endaligned
\]
We derive $l-r = f(v)$, for some $f$ one-to-one and smooth, and $w =\frac{l+r}2$. 

For classical solutions the Cauchy problem \eqref{50},\eqref{60} is equivalent to
\be
\label{80}
\left\{
\aligned
&i u_t + u_{xx} = uv + \alpha |u|^2 u
\\
& l_t - \sqrt{\sigma'(v)} l_x = (|u|^2)_x + \frac12 \Fcal(l+r)
\\
& r_t + \sqrt{\sigma'(v)} r_x = (|u|^2)_x + \frac12 \Fcal(l+r),
\endaligned
\right.
\ee
by setting $v = f^{-1}(l-r) = v(l,r)$ and with $\Fcal$ given by \eqref{70}. We take as initial data
$u(\cdot,0) = u_0 \in H^3(\RR)$, $l(\cdot,0) = l_0 \in H^2(\RR)$, $r(\cdot,0) = r_0 \in H^2(\RR)$, with
\be
\label{90}
\aligned
l_0 = w_0 + \int_0^{v_0} \sqrt{\sigma'(\xi)} \,d\xi, \qquad r_0 = w_0 - \int_0^{v_0} \sqrt{\sigma'(\xi)} \,d\xi.
\endaligned
\ee

To obtain a local strong solution of the Cauchy problem \eqref{80},\eqref{90}, we consider, using the technique employed in \cite{O} and \cite{DFO1}, an auxiliary system with non-local source terms. This is necessary in order to write the system \eqref{80},\eqref{90} without derivative loss (see \cite{DFO1}). Thus, we consider the following system,
\be
\label{100}
\left\{
\aligned
&i F_t + F_{xx} = Fv + \alpha u^2 \overline{F} + 2\alpha |u|^2 F + \frac12 u (l_x + r_x)
\\
& l_t - \sqrt{\sigma'(v)} l_x = (|\tilde u|^2)_x + \frac12 \Fcal(l+r)
\\
& r_t + \sqrt{\sigma'(v)} r_x = (|\tilde u|^2)_x + \frac12 \Fcal(l+r),
\endaligned
\right.
\ee
where $\overline{F}$ is the complex conjugate of $F$ and $F, \tilde u$ are defined by
\[
\aligned
&u(x,t) = u_0(x) + \int_0^t F(x,s) \,ds,
\\
&\tilde u(x,t) = (\Delta - 1)^{-1} (\alpha |u|^2 u + u(v-1) - i F)
\endaligned
\]
(see \cite{DFO1} for the motivation behind this definition).
The initial data are
\be
\label{110}
F(\cdot,0) = F_0 \in H^1(\RR), \qquad l(\cdot,0) = l_0 \in H^2(\RR), \qquad r(\cdot,0) = r_0 \in H^2(\RR),
\ee
with $l_0$, $r_0$ given by \eqref{90}.

We will prove the following result:
\begin{theorem}
\label{T10}
Let $(F_0, l_0, r_0) \in H^1 \times H^2 \times H^2$. Then, there exists $T^*>0$ (depending on $(F_0, l_0, r_0)$) such that for all $T<T^*$ there exists a unique solution $(F, l, r)$ of the Cauchy problem \eqref{100},\eqref{110} with 
\[
\aligned
(F, l, r) \in C^j ([0,T] ; H^{1-2j}) \times C^j ([0,T] ; H^{2-j}) \times C^j ([0,T] ; H^{2-j}), \qquad j=0,1.
\endaligned
\]
\end{theorem}

From this result and from the definitions of $F$, $\tilde u$, reasoning as in \cite[Lemma 2.1]{DFO1} (see also \cite{DFO2}), it is easy to derive the following result for the system \eqref{50},\eqref{60}:
\begin{theorem}
\label{T20}
Let $(u_0, v_0, w_0) \in  H^3 \times H^2 \times H^2$. Then, there exists $T^*>0$ (depending on $(u_0, v_0, w_0)$) such that for all $T<T^*$ there exists a unique solution $(u, v, w)$ of the Cauchy problem \eqref{50},\eqref{60} with 
\[
\aligned
(u, v, w) \in C^j ([0,T] ; H^{3-2j}) \times C^j ([0,T] ; H^{2-j}) \times C^j ([0,T] ; H^{2-j}), \qquad j=0,1.
\endaligned
\]
\end{theorem}

\proof[Proof of Theorem \ref{T10}]
In order to apply a variant of Kato's well posedness result, Theorem 6 in \cite{K}, we put the Cauchy problem \eqref{100},\eqref{110} in the framework of real spaces by introducing the new variables $F_1 = \Re F,$ $F_2 = \Im F$, $u_1 = \Re u,$ $u_2 = \Im u$. By setting $U= (F_1,F_2, l, r)$ and $F_{10} = \Re F_0,$ $F_{20} = \Im F_0$, the Cauchy problem \eqref{100},\eqref{110} can be written as follows,
\be
\label{120}
\left\{
\aligned
& U_t + A(U) U = g(t,U)
\\
& U(\cdot, 0) = U_0,
\endaligned
\right.
\ee
where
\[
A(U) = 
\left[ 
\begin{array}{cccc} 
0 & \Delta & 0 & 0 
\\
-\Delta& 0 & 0 & 0
\\ 
0 & 0 & - \sqrt{\sigma'(v)}\del_x & 0
\\
0 & 0 & 0 & \sqrt{\sigma'(v)}\del_x
\end{array} 
\right],
\]

\[
g(t,U) = 
\left[
\begin{array}{c} 
2\alpha |u|^2 F_2 - \alpha (u_1^2 - u_2^2) F_2 + 2\alpha u_1 u_2 F_1 + F_2 v + \frac12 u_2
(l_x + r_x)
\\[5pt]
2\alpha |u|^2 F_1 - \alpha (u_1^2 - u_2^2) F_1 - 2\alpha u_1 u_2 F_2 - F_1 v - \frac12 u_2
(l_x + r_x)
\\[5pt] 
(|\tilde u|^2)_x + \frac12 \Fcal(l+r)
\\[5pt]
(|\tilde u|^2)_x + \frac12 \Fcal(l+r)
\end{array} 
\right]
\]
and $U_0 = (F_{10}, F_{20} , l_0, r_0) \in Y = (H^1(\RR))^2 \times (H^2(\RR))^2$. Note that the source term is non-local. 

In what follows we use the notations of \cite[paragraph 7]{K}.
Set $X= (H^{-1}(\RR))^2 \times (L^2(\RR))^2$ and $S = (1-\Delta)\Id$, which is an isomorphism $S : Y \to X$.
Furthermore, we denote by $W_R$ the open ball in $Y$ of radius $R$ centered at the origin.

We need to check several assumptions in order to apply Theorem 6 of \cite{K}. First, it is necessary that the semigroup generated by the operator $A$ above verify 
\be
\label{125}
\aligned
\| e^{-tA(U)} \| \le e^{\omega t}, 
\endaligned
\ee 
for some real $\omega$, for all $t\ge 0$ and $U \in W_R$. Observe that it is enough to prove this only for the operator
\[
a(l,r) =
\left[ 
\begin{array}{cc} 
-\sqrt{\sigma'(v)}\del_x & 0 
\\
0 & \sqrt{\sigma'(v)}\del_x
\end{array} 
\right].
\]
since the remaining part of $A$ corresponds to the Sch\-r\"o\-din\-ger (contraction) group.
But \eqref{125} is proved in \cite[paragraph 12]{K}, with $\omega$ given by
\[
\aligned
\omega = \frac12 \sup_{x\in\RR} \| \del_x a(l,r) \| \le c(R),
\endaligned
\]
with $c : [0, +\infty) \to [0, +\infty)$ continuous.

Next, we must check that for $U \in W_R$, the property $S A(U) S^{-1} = A(U) + B(U)$ is valid for some 
$B\in \mathcal{L}(X)$. This is proved in \cite[paragraph 12]{K}: for $(l,r)$ in a ball $\widetilde W$ of $(H^2(\RR))^2$, we have
\[
\aligned
(1-\Delta) a(l,r) (1-\Delta)^{-1} = a(l,r) + B_0(l,r),
\endaligned
\]
where
\[
\aligned
B_0(l,r) = [ (1-\Delta), a(l,r)] (1-\Delta)^{-1} \in \mathcal{L}((L^2(\RR))^2)
\endaligned
\]
and
\[
B(U) = 
\left[ 
\begin{array}{cc} 
\begin{array}{cc} 
0 & 0
\\
0 & 0
\end{array} &
\begin{array}{cc} 
0 & 0
\\
0 & 0
\end{array} 
\\
\begin{array}{cc} 
0 & 0
\\
0 & 0
\end{array} &
B_0(l,r)
\end{array} 
\right],
\]
with $[\cdot,\cdot]$ denoting the matrix commutator operator.

Now, consider a pair $U,U^*\in C([0,T]; W_R),$ $U= ( F_1, F_2, l, r), U^*= ( F_1^*, F_2^*, l^*, r^*)$. It is easy to see that $g$ verifies, for fixed $T>0$, $\|g(t, U(t)\|_Y \le C(R,T),$ $t \in [0,T]$, if $ U \in C([0,T]; W_R).$ We obtain, in the same way as \cite{DFO1,DFO2},
\be
\label{130}
\aligned
\| g (\cdot, U) - g(\cdot, U^*) \|_{L^1(0,T'; X)} \le c(T') \sup_{0\le t\le T'}\| U(t) - U^*(t)\|_X
\endaligned
\ee
where $0\le T' \le T$ and $c(T')$ is a non-decreasing continuous function such that $c(0) =0$. Finally, applying Theorem~6 in \cite{K}, replacing the local condition \cite[equation (7.7)]{K} by \eqref{130}, we obtain the result. This completes the proof of Theorem~\ref{T10}.
\endproof

\section{Global existence of weak solutions for a class of stress functions}
\label{S30}
Now we will consider the question of global (in time) weak solutions of the Cauchy problem \eqref{50},\eqref{60}. For a special class of stress functions $\sigma$, we will obtain an extension of the result in \cite{NRT} for the system of nonlinear viscoelasticity and in \cite{DFF} for the system of nonlinear elasticity coupled with the nonlinear Sch\-r\"o\-din\-ger equation. We employ the adaptation of the compensated compactness method developed by D. Serre and J. Shearer in \cite{SS} for the system of nonlinear elasticity, which extends earlier results of L. Tartar \cite{T} and R.J. DiPerna \cite{DP}.

Following \cite{SS}, we define $\Sigma(v) = \int_0^v \sigma(v) \,d\xi$ (we may assume $\sigma(0) = 0$) and we impose the following conditions on the stress function $\sigma \in C^3(\RR)$:
\[
\aligned
&\text{H1: } \sigma'(v) \ge \sigma_0>0 \text{ for some constant } \sigma_0.
\\
&\text{H2: } \sigma''(\lambda_0) = 0 \text{ and } \sigma''(\lambda) \neq 0 \text{ for }\lambda\neq \lambda_0.
\\
&\text{H3: } \frac{\sigma''}{(\sigma')^{5/4}}, \quad \frac{\sigma'''}{(\sigma')^{7/4}} \in L^2(\RR),\quad
\frac{\sigma''}{(\sigma')^{3/2}}, \quad \frac{\sigma'''}{(\sigma')^{3}} \in L^\infty(\RR).
\\
&\text{H4: } \frac{\sigma(v)}{\Sigma(v)} \to 0 \text{ as } |v|\to \infty, \text{ and there exist constants $c>0$ and $q>1/2$ such that }\\& \qquad (\sigma'(v))^q \le c(1 + \Sigma(v)).
\endaligned
\]
In particular, we have $\Sigma(v) \ge \frac{\sigma_0}2 v^2$. To simplify, we also suppose that $v\in H^1(\RR)$ implies $\int_{\RR} \Sigma(v) \,dx < \infty$. A typical example is given by $\sigma(v) = v^3 +v$.

We now follow the ideas developed in \cite{DFF} and introduce a physical viscosity approximation of the Cauchy problem \eqref{50},\eqref{60}. For small $\eps >0$, we consider the system
\be
\label{200}
\left\{
\aligned
&i u_t + u_{xx} = uv + \alpha |u|^2 u
\\
&v_t = w_x && x\in\RR, \quad t\ge0,
\\
& w_t = (\sigma(v))_x + (|u|^2)_x + \Fcal(w) + \eps w_{xx},
\endaligned
\right.
\ee 
with initial data
\be
\label{210}
\aligned
u_0 ,\,   v_0 ,\,  w_0 \, \in H^1(\RR).
\endaligned
\ee
Observe that although in \cite{NRT} the authors consider initial data for which $v_0$, $w_0$ are in $L^\infty$ only, this is not possible in our case, since the method of Serre and Shearer requires initial data with stronger regularity. Of course, this does not exclude the formation of discontinuities for $t>0$.

The proof of the following lemma is an easy adaptation of the first part of the proof of Lemma~6.1 in \cite{DFF}, and so we omit it.
\begin{lemma}
Let $(u,v,w)\in \big( C([0,T); H^1(\RR) \big)^3$, $T>0$, be a solution to the Cauchy problem \eqref{200},\eqref{210} for fixed $\eps>0$. We have in $[0,T]$
\be
\label{220}
\aligned
\frac{d}{dt} \int_{\RR} |u|^2 \,dx = 0,
\endaligned
\ee
\be
\label{230}
\aligned
\frac{d}{dt} &\big[ \int_{\RR} |u_x|^2 \,dx + \frac\alpha2 \int_{\RR} |u|^4 \,dx + \int_{\RR} v|u|^2 \,dx + \frac12 \int_{\RR} w^2 \,dx + \int_{\RR} \Sigma(v) \,dx \big] 
\\
&+ \int_{\RR} \Fcal(w) w \,dx + \eps \int_{\RR} w_x^2 \,dx  = 0.
\endaligned
\ee
\end{lemma}

Now, we recall from \cite{NRT} the inequality 
\be
\label{240}
\aligned
\Big|  \int_0^t \int_{\RR} \Fcal (w) w \,dx \,d\tau \Big| \le C(T) \Big( 1 + \int_0^t \int_{\RR} w^2 \,dx \,d\tau \Big),
\endaligned
\ee
where $C(\cdot)$ is a continuous function in $[0,+\infty)$, independent of $\eps$.

We now prove some estimates which will be needed to prove the well-posedness of the viscous Cauchy problem \eqref{200},\eqref{210}.

By the Gagliardo--Nirenberg inequality	we have, by \eqref{220}, 
\[
\aligned
\int_{\RR} |u|^4 \,dx \le C_0 \| u\|^3_{L^2} \| u_x\|_{L^2} = C_0\| u_0\|^3_{L^2} \| u_x\|_{L^2},
\endaligned
\]
and since $\Sigma(v) \ge {\sigma_0} v_0^2/2$ we derive from \eqref{230},\eqref{240} for $0< t\le T$,
\be
\label{250}
\aligned
\int_{\RR} w^2(x,t) \,dx \le C(T) + C(T) \int_0^t \int_{\RR} w^2(x, \tau) \,dx \,d\tau
\endaligned
\ee
and so, by the Gronwall inequality and \eqref{230}--\eqref{250} we get for $0<t \le T$,
\be
\label{260}
\aligned
\int_{\RR} |u_x (x,t)|^2 \,dx &+ \int_{\RR} |w(x,t)|^2 \,dx + \int_{\RR} \Sigma(v(x,t)) \,dx
\\
&\qquad+ \eps \int_0^t \int_{\RR} w^2_x \,dx \,d\tau \le C(t),
\endaligned
\ee
\be
\label{270}
\aligned
\Big| \int_0^t \int_{\RR} \Fcal(w) w \,dx \,d\tau \Big| \le C(T),
\endaligned
\ee
for some continuous function $C$ on $[0, +\infty)$ independent of $\eps$.

Next, we will prove the following estimate, for $\eps\le 1$, $t \in [0,T]$:
\be
\label{280}
\aligned
\eps \int_0^t \int_{\RR} \sigma'(v) (v_x)^2 \,dx \,d\tau + \eps^2 \int_{\RR} (v_x)^2 \,dx \le C(T),
\endaligned
\ee
where $C$ is a continuous function on $[0, +\infty)$ independent of $\eps$. For this purpose we follow the ideas in \cite[proof of (8)]{SS} (see also \cite{DFF}). We deduce from \eqref{200}, for $v$ smooth enough,
\[
\aligned
\int_{\RR} w_t v_x - \sigma'(v) (v_x)^2 \,dx = \int_{\RR} (|u|^2)_x v_x \,dx + \int_{\RR} \Fcal(w) v_x \,dx +\eps \int_{\RR} w_{xx} v_x \,dx
\endaligned
\]
and 
\[
\aligned
- \frac d{dt} \int_{\RR} w_x v \,dx &+ \int_{\RR} (w_x)^2 \,dx - \int_{\RR} \sigma'(v) (v_x)^2 \,dx 
\\
&= \int_{\RR}(|u|^2)_x v_x \,dx + \int_{\RR} \Fcal(w) v_x + \frac\eps2 \frac d{dt} \int_{\RR} (v_x)^2 \,dx
\endaligned
\]
(since, for $w$ smooth enough, $-\frac d{dt} \int_{\RR} w_x v \,dx = -  \int_{\RR} w_{xt} v \,dx - \int_{\RR} w_x v_t \,dx = \int_{\RR} w_t v_x \,dx - \int_{\RR} w_x v_t \,dx $ and $v_t = w_x$).

Integrating over $(0,t)$ we obtain, with $v_0(x) = v(x,0),$ $w_0(x) = w(x,0)$,
\[
\aligned
- \int_{\RR}& w_x v \,dx + \int_{\RR} w_{0x} v_0 \,dx + \int_0^t \int_{\RR} (w_x)^2 \,dx \,d\tau - \int_0^t 
\int_{\RR} \sigma'(v) (v_x)^2 \,dx \,d\tau 
\\
&= \int_0^t \int_{\RR} (|u|^2)_x v_x \,dx \,d\tau 
+\int_0^t \int_{\RR}  \Fcal(w) v_x \,dx \,d\tau + \frac\eps2 \int_{\RR} (v_x)^2 - (v_{0x})^2 \,dx.
\endaligned
\]
From $ -\int_{\RR} w_x v \,dx = \int_{\RR} w v_x$ we derive
\be
\label{290}
\aligned
&\int_0^t \int_{\RR} \sigma'(v) (v_x)^2 \,dx \,d\tau + \frac\eps2 \int_{\RR} (v_x)^2 \,dx 
\le \frac\eps4
\int (v_x)^2 \,dx + \frac1\eps \int_{\RR} w^2 \,dx 
\\
&\quad+ \int_{\RR} |w_{0x}| |v_0| \,dx + \frac\eps2 
\int_{\RR} (v_{0x})^2 \,dx 
+ \eps \int_0^t \int_{\RR}  (w_x)^2 \,dx \,d\tau 
\\
&\quad+ 2 \int_0^t \int_{\RR}  |u| |u_x| |v_x| \,dx \,d\tau + \Big| \int_0^t \int_{\RR}  \Fcal(w) v_x \,dx \,d\tau \Big|.
\endaligned
\ee
Moreover, we easily get from \eqref{260}, for a fixed $\delta>0$,
\be
\label{300}
\aligned
2 \int_0^t \int_{\RR}  |u| |u_x| |v_x| \,dx \,d\tau \le C(\delta) C(T) +  \frac\delta{\sigma_0} 
\int_0^t \int_{\RR}  \sigma'(v) (v_x)^2 \,dx \,d\tau,
\endaligned
\ee
and also by the definition of $\Fcal$,
\be
\label{310}
\aligned
\Big| \int_0^t \int_{\RR} & \Fcal(w) v_x \,dx \,d\tau \Big| 
\le \delta \int_0^t \int_{\RR}  (v_x)^2 \,dx \,d\tau + C(\delta) C(T) \int_0^t \int_{\RR} w^2 \,dx \,d\tau
\\
&\le \frac\delta{\sigma_0} 
\int_0^t \int_{\RR}  \sigma'(v) (v_x)^2 \,dx \,d\tau + C(\delta) C(T) \int_0^t \int_{\RR} w^2 \,dx \,d\tau.
\endaligned
\ee
From \eqref{290}--\eqref{310} we derive, choosing $\delta = \frac{\sigma_0}4$ and multiplying by $\eps$, the estimate \eqref{280} for $t\in (0,T]$.

Let us now analyze the problem of the existence of a unique solution $(u,v,w) \in (C ([0,+\infty); H^1))^3$ of the Cauchy problem \eqref{200},\eqref{210} for fixed $\eps \in(0,1]$. We assume $\eps =1$ without loss of generality. Let us begin with the study of the existence and uniqueness of a local (in time) solution. Let $0<T<+\infty$ and introduce the complex Banach space (not to be confused with the spaces introduced in Section~\ref{S20}) $X_T = C([0,T]; H^1)$ and the real Banach space $Y_T = C([0,T]; H^1)$, with the standard norms. Let us consider the product space $\BB_R^T \times (B_R^T)^2$ where 
$\BB_R^T = \{ u \in X_T : \| u\|_{X_T} \le R \}$ and $(B_R^T)^2 = \{ (v,w) \in Y_T^2 : \| v\|_{Y_T} \le R, \| w\|_{Y_T} \le R\}$.

Given $(\tilde u, \tilde v, \tilde w) \in \BB_R^T \times (B_R^T)^2$ we consider the mapping
\[
\aligned
(\tilde u, \tilde v, \tilde w) \mapsto \Phi (\tilde u, \tilde v, \tilde w) = (u, v, w)\in X_T \times (Y_T)^2
\endaligned
\]
where $(u, v, w)$ is the unique solution of the linear problem
\be
\label{320}
\left\{
\aligned
&i u_t + u_{xx} = \alpha | \tilde u|^2 \tilde u + \tilde v \tilde u, &&u(0) = u_0,
\\
&w_t - w_{xx} = (\sigma(\tilde v))_x + (|\tilde u|^2)_x + \Fcal(\tilde w), && w(0) = w_0,
\\
& v(t) = v_0 +\int_0^t w_x \,d\tau, 
\endaligned
\right.
\ee
with $u_0, v_0, w_0 \in H^1$ and such that $w \in L^2(0,T; H^2), w_t \in L^2(0,T;L^2)$.
With the usual method of semigroups, we have
\[
\aligned
u(t) = e^{i\Delta t} u_0 - i \int_0^t e^{i \Delta (t-s)} (\alpha |\tilde u|^2 \tilde u + \tilde v\tilde u)(s) \,ds
\endaligned
\]
and 
\[
(\sigma(\tilde v))_x + (|\tilde u|^2)_x + \Fcal(\tilde w) \in C([0,T]; L^2)
\]
with the estimate \cite[(3.42)]{NRT}
\be
\label{330}
\aligned
\int_0^t \int_{\RR} \big| \Fcal (\tilde w_1) - \Fcal (\tilde w_2) \big|^2 \,dx \,d\tau\le C(T) \int_0^t \int_{\RR} | \tilde w_1 - \tilde w_2 |^2 \,dx \,d\tau.
\endaligned
\ee
Reasoning as in \cite[proof of Theorem~3.1]{NRT}, we can, for convenient $T$ and $R > \max (|u_0|_{H^1}, |v_0|_{H^1}, |w_0|_{H^1})$, apply the Banach fixed point theorem for strict contractions to obtain existence and uniqueness of a local (in time) solution of the viscous system \eqref{200},\eqref{210}. Moreover, we have
\[
\aligned
w_t - w_{xx} = (\sigma(v))_x + (|u|^2)_x + \Fcal(w)
\endaligned
\]
and, from the previous estimates \eqref{260},\eqref{280}, \eqref{330} (with $\tilde w_1 = w,$ $\tilde w_2 =0$) and \eqref{220}, we derive the \emph{a priori} estimate
\[
\aligned
\|w_t - w_{xx} \|_{L^2(0,T;L^2)} \le c(T),
\endaligned
\]
with $c \in C([0,+\infty[ ; \RR_+)$, and this implies, by standard regularity estimates for parabolic equations, that $w \in L^2(0,T; H^2)$ and $\|w \|_{L^2(0,T; H^2)} \le c(T)$. This in turn gives similar 
\emph{a priori} estimates for $\|w_t \|_{L^2(0,T; L^2)}$ and $\|w \|_{C([0,T]; H^1)}.$ 

We are now in a position to pass to the limit $\eps\to 0$ and state the main result in this section.

\begin{theorem}
\label{T30}
Assume $(u_0, v_0, w_0) \in (H^1(\RR))^3$ and let $L^\Sigma_\loc(\RR \times [0,+\infty)) $ be the space of functions $v$ measurable in $\RR \times [0,+\infty)$ such that
\[
\aligned
\int_K \Sigma(v) \,dx \,dt < \infty \text{ for each compact } K \subset \RR\times [0,+\infty)\endaligned.
\]
Then there exists 
\[
\aligned
(u,v,w) \in L^\infty_\loc ((0, +\infty); H^1) \times L^\Sigma_\loc (\RR \times [0, +\infty)) \times L^2_\loc (\RR \times [0,+\infty))
\endaligned
\]
such that 
\[
\aligned
- i & \int_0^\infty \int_{\RR} u \theta_t \,dx \,dt - \int_0^\infty \int_{\RR} u_x \theta_x \,dx \,dt + \int_{\RR} u_0(x) \theta(x, 0) \,dx 
\\
&= \int_0^\infty \int_{\RR} v u \theta \,dx \,dt + \alpha \int_0^\infty \int_{\RR} |u|^2 u \theta \,dx \,dt,
\endaligned
\]
for all complex-valued $\theta \in C^1_0 (\RR\times [0, +\infty))$, and
\[
\aligned
&\int_0^\infty \int_{\RR} (v \phi_t - w \phi_x) \,dx \,dt  + \int_{\RR} v_0(x) \phi(x,0) \,dx 
\\
&\quad+ \int_0^\infty \int_{\RR} w \psi_t - \sigma(v) \psi_x + \Fcal(w) \psi \,dx \,dt 
\\
&\qquad+ \int_{\RR} w_0(x) \psi(x,0) \,dx + \int_0^\infty \int_{\RR} (|u|^2)_x \psi \,dx \,dt =0,
\endaligned
\]
for every real-valued $\phi, \psi \in C^1_0 (\RR\times [0, +\infty))$.
\end{theorem}

\proof
We follow the ideas in \cite{SS,DFF} (see also \cite{C}): for each $\eps \in (0,1]$ let $(u_\eps, v_\eps, w_\eps)\in \big( C([0,+\infty);H^1)\big)^3$ be the unique solution of the Cauchy problem \eqref{200},\eqref{210}. From the system \eqref{200}, \eqref{220}, and \eqref{260}, we derive, for fixed $T>0$,
\[
\aligned
& \{ u_\eps \}_\eps \text{ bounded in } L^\infty ((0,T); H^1),
\\
& \{ u_{\eps t} \}_\eps \text{ bounded in } L^\infty ((0,T); H^{-1}),
\endaligned
\]
and so, by a well known compactness result, $\{ u_\eps \} $ is in a compact set of 
$L^2( 0,T; L^2(I_R))$ for each interval $I_R = (-R,R)$, $R>0$. Furthermore, there exists $u\in H^1(\RR)$ such that (for a subsequence) $u_\eps \ws u$ in $L^\infty_\loc(0,+\infty; H^1(\RR))$ when $\eps \to 0$.
By a standard diagonalization method, we conclude that in fact, for a suitable subsequence, 
$u_\eps \to u$ in $L^1_\loc(\RR \times [0, +\infty))$.

We also have, by \eqref{260}, 
\[
\aligned
& \{ w_\eps \}_\eps \text{ bounded in } L^2_\loc(\RR \times [0, +\infty)),
\\
& \{ v_{\eps} \}_\eps \text{ bounded in } L^\Sigma(\RR \times [0, +\infty)).
\endaligned
\]
By \eqref{260} and \eqref{280} we derive
\be
\label{340}
\aligned
\eps \int_0^t \int_{\RR} (w_{\eps x} )^2 + \sigma'(v_\eps) (v_{\eps x})^2 \,dx \,dt \le C(T)
\endaligned
\ee
where $C(t)$ is a continuous function depending the $H^1$ norm of the initial data, but not on $\eps$.

Now consider the quasilinear system
\be
\label{350}
\left\{
\aligned
& v_t = w_x
\\
&w_t = (\sigma(v))_x
\endaligned\right.
\ee
and, for $(v,w) \in \RR^2$, let $\eta(v,w), q(v,w)$ be a pair of smooth convex entropy-entropy flux pair for \eqref{350} such that $\eta_w$, $\eta_{ww}$ and $\eta_{vw}/\sqrt{\sigma}$ are bounded functions in $\RR^2$. From \eqref{220},\eqref{260},\eqref{340} and since, for each interval $I_R=[-R,R]$, $R>0$, we have
\[
\aligned
\Big| \int_0^t \int_{I_R} \Fcal(w_\eps) \,dx \,dt \Big| \le C_R(t),
\endaligned
\]
where $C_R(t)$ is a continuous function on $[0, +\infty)$ independent of $\eps$, we derive from the system \eqref{200} (cf.~similar estimates in \cite{SS,C} and \cite{DFF}) that
\[
\aligned
\del_t \eta(v_\eps, w_\eps) + \del_x q(v_\eps, w_\eps)
\endaligned
\]
belongs to a compact set of $W^{-1,2}_\loc (\RR\times [0,+\infty)).$ Hence, in view of the assumptions H1--H4, we may apply the result on compensated compactness of D.~Serre and J.~Shearer \cite{SS} to conclude that $\{ (v_\eps, w_\eps) \}_\eps$ is pre-compact in $\big( L^1_\loc(\RR \times [0, +\infty)) \big)^2$. Hence, there exist a subsequence $\{ (u_\eps, v_\eps, w_\eps) \}_\eps$ and 
\[
\aligned
(u,v,w) \in L^\infty_\loc(0,+\infty; H^1) \times L^\Sigma_\loc(\RR \times [0, +\infty)) \times L^2_\loc(\RR \times [0, +\infty))
\endaligned
\]
such that 
\[
\aligned
(u_\eps, v_\eps, w_\eps) \to (u,v,w) \text{ in  }\big( L^1_\loc(\RR \times [0, +\infty)) \big)^3, \quad \eps\to 0.
\endaligned
\]
Now take functions $\theta, \phi$ and $\psi$ as introduced in the statement of Theorem~\ref{T30}, multiply the first equation of system \eqref{200} by $\theta$, the second one by $\phi$ and the third one by $\psi$ and integrate over $\RR \times [0, +\infty)$. After integrating by parts, we pass to the limit as $\eps \to 0$ in view of the above convergences. This allows us to prove that $(u,v,w)$ are weak solutions of the system \eqref{10}, since we have for some $T= T(\psi)$,
\[
\aligned
\Big| \eps \int_0^\infty \int_{\RR} w_{\eps xx} \psi \,dx \,dt \Big| \le c_\psi \eps^{1/2} \Big( \eps \int_0^T \int_{\RR} (w_{\eps x})^2 \,dx \,dt \Big)^{1/2},
\endaligned
\]
which goes to zero as $\eps \to 0$ by \eqref{260}. This completes the proof of Theorem~\ref{T30}.
\endproof

\section{Numerical experiments}
\label{S40}

In this section we present some numerical experiments on the system \eqref{10}, to illustrate our results. 
We use the more convenient formulation \eqref{50}--\eqref{70}.
According to \cite{NRT}, if the kernel $k$ in \eqref{10} is the derivative of a $C^2$, positive, decreasing and convex function on $[0,+\infty)$, then the term $q(0)w(x,t)$ in $\Fcal$ has a damping effect and the solutions should remain classical and globally defined.

Thus, we will consider the particular case $k(t) = e^{-t}$ in \eqref{10}, which gives, after some elementary calculations, $q(t) = e^{-2t}$ in \eqref{50}.

The stress function is taken as $\sigma(v) = v^3 + v$, which verifies the conditions H1--H4 in Section~\ref{S30}. The memory term in \eqref{50} reads
\be
\label{390}
\aligned
\Fcal(w) = w(x,t) - e^{-2t} w_0(x) - 2\int_0^t e^{2(s-t)} w(x,s) \,ds.
\endaligned
\ee

\subsection{The numerical scheme}
We use a fourth-order explicit Runge--Kutta scheme for the time-stepping along with standard finite difference discretizations of the space derivatives. As is usual in the simulation of problems posed on the whole line, we restrict ourselves to a bounded domain $[L_1, L_2]$ and to initial data decaying exponentially for large $|x|$, and set the boundary conditions to zero. 

Thus, we are given a spatial mesh size $h >0$, $J = (L_2 - L_1)/h$ (which we suppose is a natural number without loss of generality), a time step $\tau$, a suitable discretization of the initial data, $(u_{0j}, v_{0j}, w_{0j})_{j=0,\dots,J},$ an approximation $(\unj, \vnj, \wnj)_{j=0,\dots,J}$ of $(u(x_j, t_n), v(x_j, t_n), w(x_j, t_n))$ at time $t_n = n\tau,$ $ n\in \NN$, and $x_j = L_1 + jh $.  We then obtain $(\unnj, \vnnj, \wnnj)$ by solving by the fourth-order Runge--Kutta method the equations 
\be
\label{400}
\aligned
&\frac d{dt} \unj = \frac i{h^2} (\unjj - 2 \unj + \unjm) - i \unj |\unj|^2 - i \vnj \unj
\\
&\frac d{dt} \vnj = \frac1{2h}(\wnjj - \wnjm)
\\
&\frac d{dt} \wnj = \frac1{2h}( \sigma(\vnjj) - \sigma(\vnjm)) + \frac{h}2 (\wnjj - 2\wnj + \wnjm) 
\\
&\qquad+ \frac1{h^2} ( |\unjj|^2 - |\unjm|^2) + \wnj - e^{-2t_n} w_{0j} -  2 e^{-2t_n} \Fcal_j^n,
\endaligned
\ee
where $\Fcal_j^n$, the discretization of the non-local integral term in \eqref{390} 
is defined recursively by
\[
\aligned
\Fcal_j^n =  \Fcal_j^{n-1} + \tau (e^{2t_{n-1}} w_j^{n-1}).
\endaligned
\]
Note the introduction of a viscosity term in the third equation on \eqref{400} to improve stability, as is customary in the numerical treatment of hyperbolic equations. 

As initial data (see Figure~\ref{420}), we have set
\be
\label{410}
\aligned
&u_0 (x) = C e^{30 i x}  \cosh^{-1} (\sqrt{50} x),
\\
&v_0 (x) = C\cosh ^{-1}(\sqrt{20} (x-0.1)),
\\
&w_0 (x) = C\cosh ^{-1}(\sqrt{20} (x+0.1)).
\endaligned
\ee
For other simulations of short wave long wave interactions in the case of a Sch\-r\"o\-din\-ger equation coupled with a nonlinear conservation law, see \cite{AF1,AF2}.

\subsection{Numerical results}
\begin{figure}
\includegraphics[width=\linewidth,keepaspectratio=false]{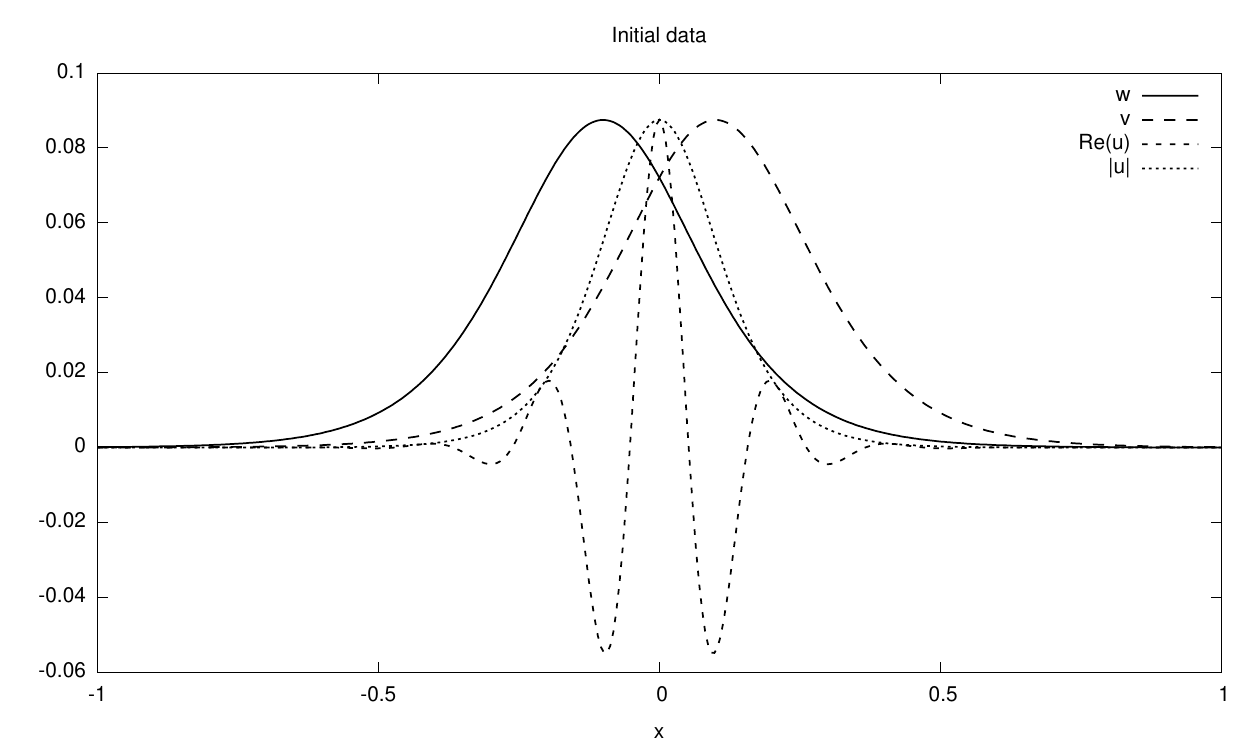}
\caption{The initial data \eqref{410}.}
\label{420}
\end{figure}

\begin{figure}
\includegraphics[width=\linewidth,keepaspectratio=false]{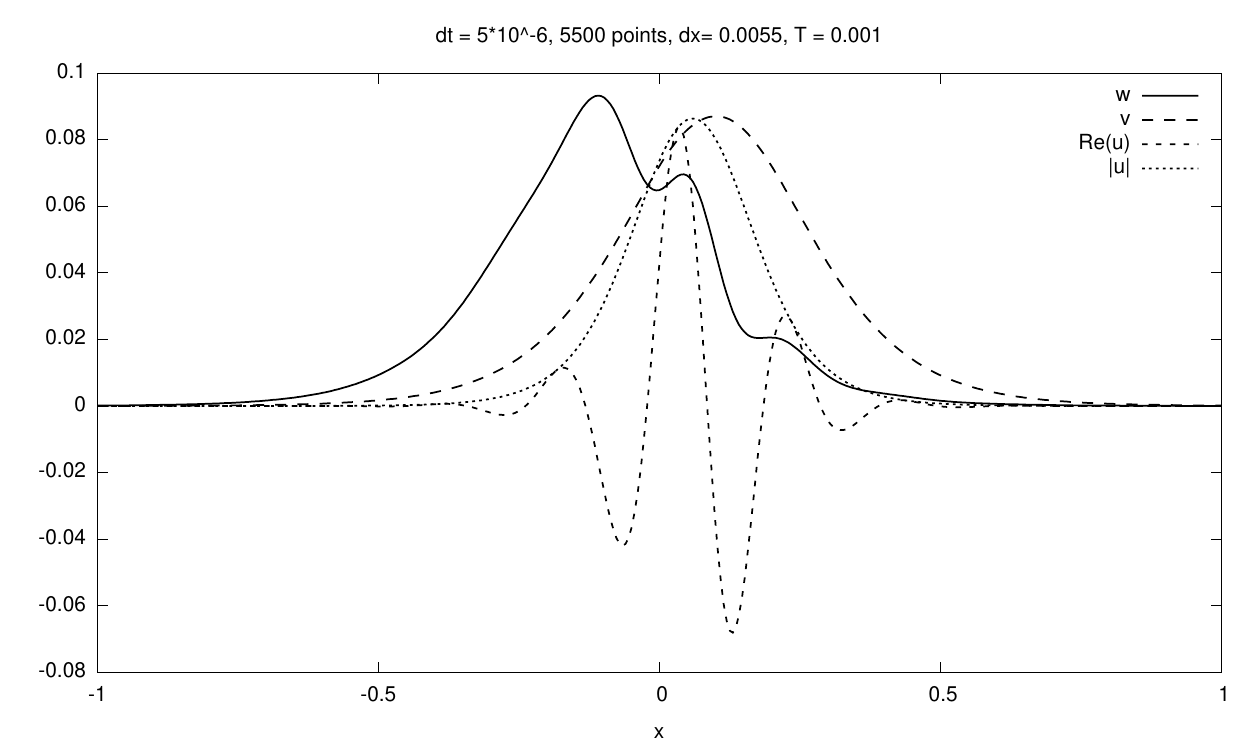}
\caption{Computed solution of system \eqref{50},\eqref{60}, $T=0.001$.}
\label{430}
\end{figure}

\begin{figure}
\includegraphics[width=\linewidth,keepaspectratio=false]{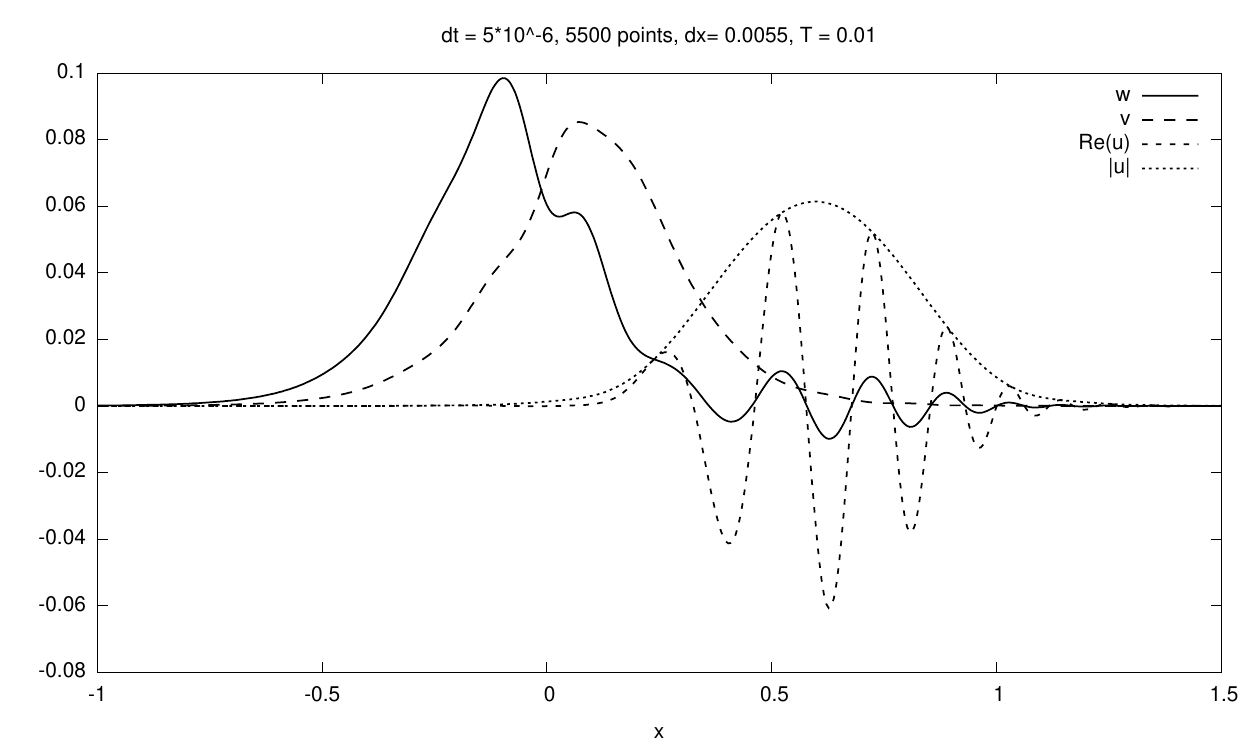}
\caption{Computed solution of system \eqref{50},\eqref{60}, $T=0.01$.}
\label{440}
\end{figure}

\begin{figure}
\includegraphics[width=\linewidth,keepaspectratio=false]{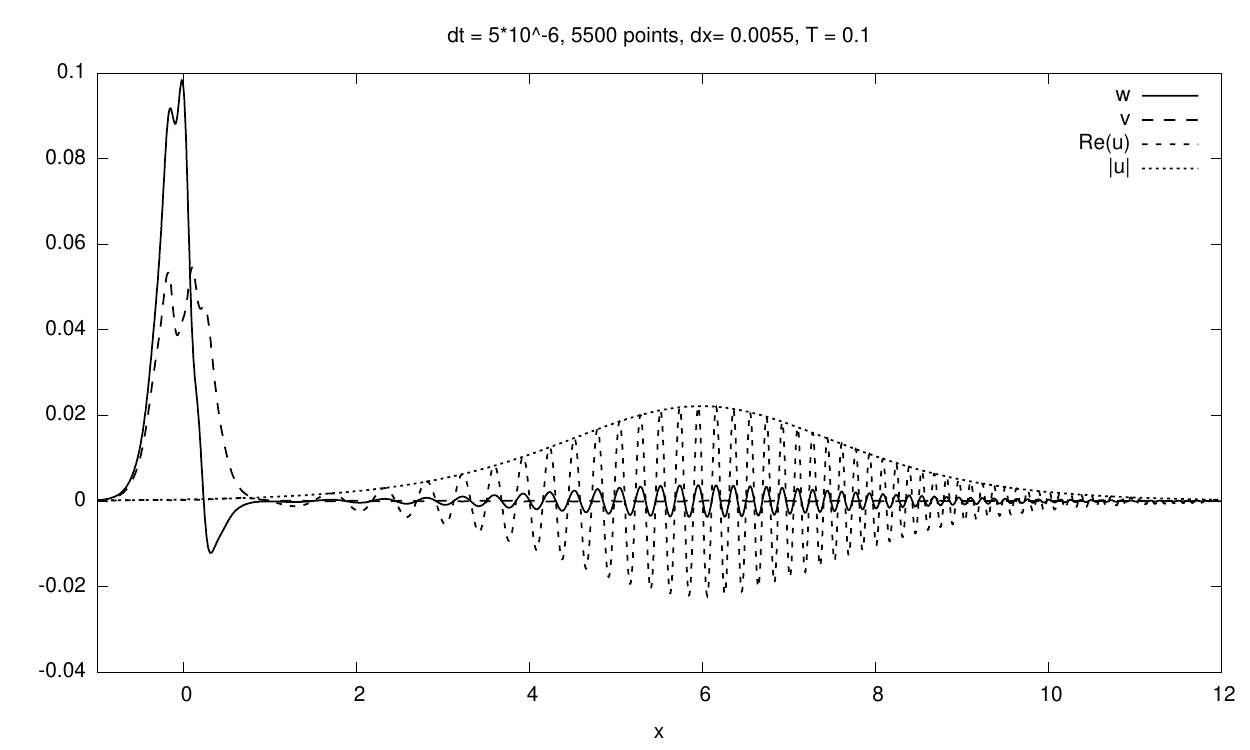}
\caption{Computed solution of system \eqref{50},\eqref{60}, $T=0.1$.}
\label{450}
\end{figure}

\begin{figure}
\includegraphics[width=\linewidth,keepaspectratio=false]{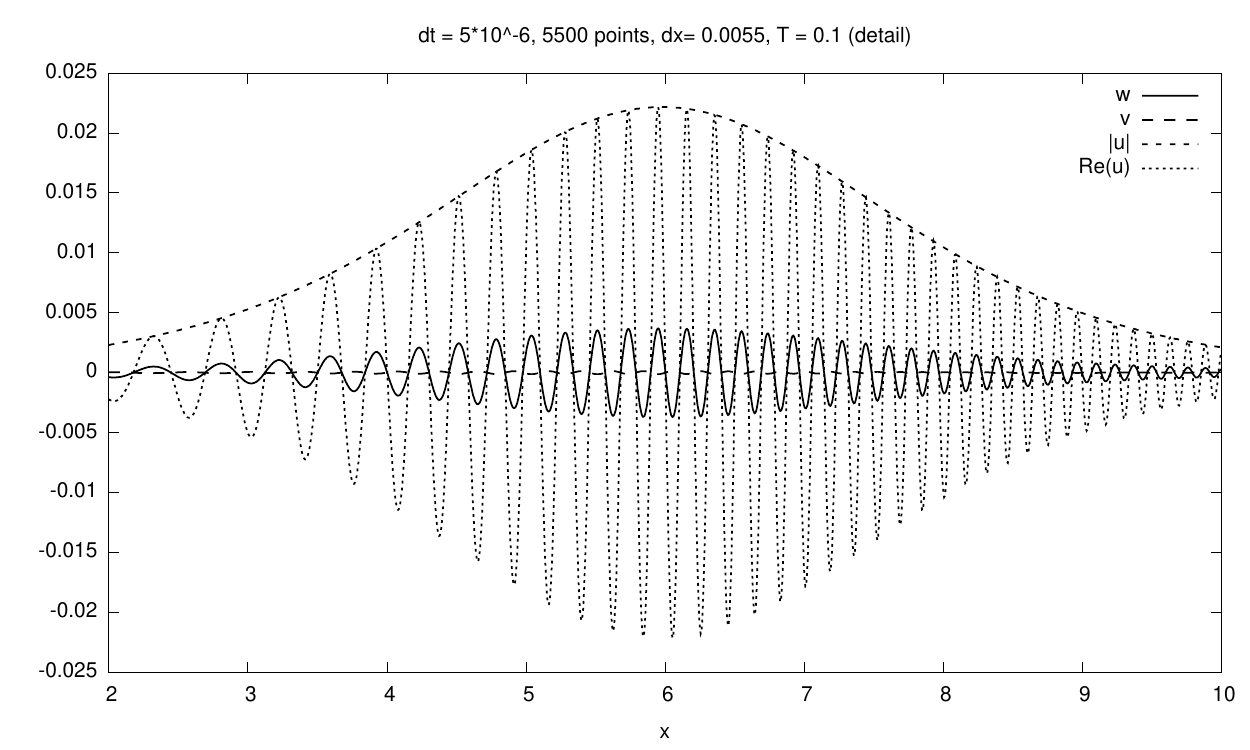}
\caption{Computed solution of system \eqref{50},\eqref{60}, $T = 0.1$ (close up).}
\label{460}
\end{figure}

We now present some results obtained by this numerical method. In Figure~\ref{420}, we present the initial data \eqref{410}. Next, In Figures~\ref{430}--\ref{460}, we plot the numerical solution, respectively, at times $T = 0.001,$ $T = 0.01,$ and $T = 0.1$, with a close-up of the interaction region in 
Figure~\ref{460}.

We can clearly observe in these simulations the interaction between the short and the long waves. In particular, there is the formation of a train of small frequency waves in the $v$ and (more strongly) in the $w$ variables which can clearly be discerned. These new waves, due to the nonlinear coupling, are coherent with the oscillations in the Sch\-r\"o\-din\-ger variable $u$, as can be seen in Figures~\ref{440}--\ref{460}. 

Thus, our simulations show creation of new waves due to the interaction effects and thus allow us to see a new qualitative property of the solutions to the system under consideration.

\section*{Acknowledgements}
The authors would like to thank Prof. Jos\'e Teixeira of Centre de l'\'En\'ergie Atomique, France, for helpful comments.
The authors 
were supported by FCT, through Financiamento Base 2008-ISFL-1-209 and the FCT grant PTDC/MAT/110613/2009. PA was also supported by 
FCT through a {Ci\^encia 2008} fellowship.


\end{document}